\documentclass[12pt,a4paper]{article}

\usepackage{makeidx}
\usepackage{amssymb}
\usepackage{amsfonts}
\usepackage{amsmath}
\usepackage{graphicx}

\newtheorem{theorem}{Theorem}[section]

\newtheorem{condition}{Condition}

\newtheorem{definition}{Definition}[section]

\newtheorem{remark}{Remark}[section]

\title{The volume and time comparison principle and transition probability
estimates for random walks}

\author{Andr\'{a}s Telcs \\
{ Department of Computer Science and Information Theory, }\\
{ Budapest} { University of Technology and Economics}\\
{ telcs@szit.bme.hu}}

\begin{document}

\maketitle

\begin{abstract}
This paper presents necessary and sufficient conditions for on-
and off-diagonal \ transition probability estimates for random
walks on weighted graphs. \ On the integer lattice and on may
fractal type graphs both the volume of a ball and the mean exit
time from a ball is independent of the centre, uniform in space.
Here the upper estimate is given without such restriction and
two-sided estimate is given if uniformity in the space assumed
only for the mean exit time.
\end{abstract}
\section{ Introduction}

\setcounter{equation}{0}\label{sintr}

This paper presents on- and off-diagonal transition probability estimates on
weighted graphs. \ The central object of the investigation is the minimal
solution of the discrete heat equation
\begin{equation}
\triangle u_{n}=\frac{\partial }{\partial n}u_{n}
\end{equation}
on weighted graphs, where $\triangle =P-I$ is the discrete Laplace and $%
\frac{\partial }{\partial n}u_{n}=u_{n+1}-u_{n}$ is the differential
operator. \ The classical form of the diagonal upper estimate for the
minimal solution is
\begin{equation}
p_{n}\left( x,x\right) \leq Cn^{-\frac{d}{2}},
\end{equation}
which holds on $Z^{d},d \in N$ , where $V\left(
x,R\right) $ the volume of a ball of radius $R$ has uniformly
polynomial growth $V\left( x,R\right) \simeq R^{d}.$ \ Coulhon
and Grigor'yan \cite{CG} proved diagonal estimates for graphs for
non-uniform volume. They have shown that
\begin{equation}
p_{n}\left( x,x\right) \leq \frac{C}{V\left( x,\sqrt{n}\right) }
\textnormal{ and the volume doubling condition}
\end{equation} is equivalent to
\begin{equation}
p_{n}\left( x,y\right) \leq \frac{C}{V\left( x,\sqrt{n}\right) }%
exp\left( -c\frac{d^{2}\left( x,y\right) }{n}\right)
\textnormal{and the volume doubling condition}  \label{cg}
\end{equation}%
and is equivalent to the relative Faber-Krahn inequality. For a
detailed introduction and history see Woess \cite{W} (or Barlow
\cite{B}, Coulhon \cite{Co1} , Grigor'yan \cite{Gr1}, Varopoulos,
Saloff-Coste, Coulhon \cite {VCS}).

\ Another branch of the generalization was the study of fractals and fractal
like graphs cf. \cite{B} to obtain sub-Gaussian estimates
\begin{equation}
p_{n}\left( x,y\right) \simeq \frac{C}{n^{\frac{\alpha }{\beta }}}exp\left( -%
\frac{d^{\beta }\left( x,y\right) }{Cn}\right) ^{\frac{1}{\beta -1}}
\end{equation}
where the constant $C$ different for the upper and lower bound. \ Here the
sub-Gaussian feature is provided by the $\beta >2$ exponent which describes
the mean exit time
\begin{equation}
E\left( x,R\right) \simeq R^{\beta }
\end{equation}
of a ball $B\left( x,R\right) .$ \ Let us consider the (generalized) inverse
function $e\left( x,n\right) $ of the mean exit time $E\left( x,R\right) $
in the second variable. \ Based on the usual heuristic one might expect that
a fully local diagonal upper estimate of the form of
\begin{equation}
p_{n}\left( x,x\right) \leq \frac{C}{V\left( x,e\left( x,n\right) \right) }
\end{equation}
can be given. This paper announces on- and off- diagonal estimates of this
local type. The sub-Gaussian exponents of the off-diagonal upper and lower
estimates do not coincide in this generality ( for further explanation and
examples see \cite{HK}). \ It can be seen that the sub-Gaussian exponents
meet if and only if the mean exit time is uniform in the space i.e.
\begin{equation}
E\left( x,R\right) \simeq F\left( R\right)
\end{equation}
for a function $F.$ This is usually called in the physics
literature space-time scale function and a semi-local framework
can be developed in its presence.

\section{Preliminaries}

Let us consider a countable infinite connected graph $\Gamma $. \
A weight function $\mu _{x,y}=\mu _{y,x}> 0$ is given on the
edges $x\sim y.$ This
weight induces a measure $\mu (x)$%
\begin{equation}
\mu (x)=\sum_{y\sim x}\mu _{x,y}
\end{equation}
\begin{equation}
\mu(A)=\sum_{y\in A}\mu (y)
\end{equation}
on the vertex set $A\subset \Gamma $ and defines a reversible Markov chain $%
X_{n}\in \Gamma $, i.e. a random walk on the weighted graph $(\Gamma ,\mu )$
with transition probabilities
\begin{equation}
P(x,y) =\frac{\mu_{x,y}}{\mu_(x)},
\end{equation}
\begin{equation}
P_{n}(x,y) =P(X_{n}=y|X_{0}=x).
\end{equation}

\begin{condition}
In the whole sequel we assume that condition
$\mathbf{(p}_{0}\mathbf{)} $ holds, that is, there is a universal
$p_{0}>0$ such that for all $x,y\in \Gamma ,x\sim y$
\begin{equation}
\frac{\mu _{x,y}}{\mu (x)}\geq p_{0}.  \label{p0}
\end{equation}
\end{condition}

The graph is equipped with the usual (shortest path length) graph distance $%
d(x,y)$ and open metric balls are defined for $x\in \Gamma ,$ $R>0$ as $%
B(x,R)=\{y\in \Gamma :d(x,y)<R\}$ and its $\mu -$measure is denoted by $%
V(x,R)$.

\begin{definition}
The weighted graphs satisfies the volume comparison principle %
$\left( \mathbf{VC} \right) $ (c.f. \cite{GM}) if there is a
constant $C_{V}>1$
such that for all $x\in \Gamma $ and $R>0,y\in B\left( x,R\right) $%
\begin{equation}
\frac{V(x,2R)}{V\left( y,R\right) }\leq C_{V}.  \label{VC}
\end{equation}
\end{definition}

\begin{definition}
The weighted graph has the \textbf{volume doubling} $(\mathbf{VD)}$
 property if there is a constant $%
D_{V}>0$ such that for all $x\in \Gamma $ and $R>0$%
\begin{equation}
V(x,2R)\leq D_{V}V(x,R).  \label{PDV1}
\end{equation}
\end{definition}

One can see that $\left( VD\right) $\ and $\left( VC\right) $ are equivalent.%
\newline

\section{Upper estimates}

This section is mainly devoted to upper estimates, but at the end lower
estimates are also given \ providing comparison with the upper one..

Let us consider the exit time $T_{B(x,R)}=\min \{k:X_{k}\notin B(x,R)\}$
from the ball $B(x,R)$ and its mean value $E_{z}(x,R)=E%
(T_{B(x,R)}|X_{0}=z)$ and let us use the $E(x,R)=E_{x}(x,R)$ short notation.
In the analogy to the volume comparison we introduce the (mean exit) time
comparison principle.

\begin{definition}
We will say that the weighted graph $(\Gamma ,\mu )$ satisfies the time comparison principle $\left( \mathbf{TC}\right) $ if there is a
constant $C_{T}>1$ such that for all $x\in \Gamma $ and $R>0,y\in B\left(
x,R\right) $%
\begin{equation}
\frac{E(x,2R)}{E\left( y,R\right) }\leq C_{T}.  \label{TC}
\end{equation}
\end{definition}

\begin{definition}
We will say that $(\Gamma ,\mu )$ has the \textbf{time doubling} property %
$\left( \mathbf{TD}\right) $ if there is a $D_{T}>0$
such that for all $x\in \Gamma $ and $R\geq 0$%
\begin{equation}
E(x,2R)\leq D_{T}E(x,R).  \label{TD}
\end{equation}
\end{definition}

One should notice that $\left( TC\right) $ implies $\left(
TD\right)  $ but the opposite is not true in general.

Basically the $\left( VC\right) $ and $\left( TC\right) $
principles specify the framework of the \textit{local setup }for
our study.

We introduce the skewed version of the parabolic mean value inequality.

\begin{definition}
We shall say that the skewed parabolic mean value inequality $\left(
sPMV\right) $ holds if there are $%
0<c_{_{1}}<c_{_{2}}\leq 1\leq C$ such that for all $R>0,x\in \Gamma ,y\in
B\left( x,R\right) $ for all non-negative solutions $u_{n}$ of the discrete
heat equation on $\left[ 0,c_{2}E\left( x,R\right) \right] \times B\left(
x,R\right) $%
\begin{equation}
u_{n}(x)\leq \frac{C}{V(y,2R)E\left( y,2R\right) }\sum_{i=c_{1}E}^{n}\sum_{z%
\in B(x,R)}u_{i}(z)\mu (z)\,  \label{sPMV}
\end{equation}
satisfied, where $E=E\left( x,R\right) ,n=c_{2}E$.
\end{definition}

\begin{definition}
We shall say that the mean value inequality $\left( \mathbf{MV}\right)
 $ holds, i.e. for all $x\in \Gamma ,R>0$ and for
all function $u\geq 0$ on $\overline{B}\left( x,R\right) $ which is harmonic
on $B=B\left( x,R\right) $%
\begin{equation}
u\left( x\right) \leq \frac{C}{V\left( x,R\right) }\sum_{y\in B}u\left(
y\right) \mu \left( y\right) .  \label{mv}
\end{equation}
\end{definition}

\begin{definition}
The local kernel function $1\leq k=$ $k_{y}=k_{y}(n,R)\leq n,$ is defined as
the maximal integer for which
\begin{equation}
\frac{n}{k}\leq qE(y,\left\lfloor \frac{R}{k}\right\rfloor )
\end{equation}
or $k=0$ by definition if there is no appropriate $k$. Here $q$ is a small
fixed constant.
\end{definition}

\begin{definition}
For convenience we will use the following notation
\begin{equation}
k_{C}\left( x,n,R\right) =\min\limits_{z\in B\left( x,e\left( x,n\right)
\right) }\left\{ k_{z}\left( Cn,\frac{1}{C}R\right) \right\} .
\end{equation}
Denote $R=d\left( x,y\right) $ and let
\begin{equation}
\kappa _{C}\left( n,x,y\right) =\max \left\{ k_{C}\left( x,n,R\right)
,k_{C}\left( y,n,R\right) \right\}  \label{kappdef}
\end{equation}
if $r>3\left[ e\left( x,n\right) +e\left( y,n\right) \right] $ \ and $\kappa
_{C}=0$ otherwise.
\end{definition}

\begin{theorem}
\label{tLUE} If $\left( \Gamma ,\mu \right) $ satisfies $\left( p_{0}\right)
$ then the following conditions are equivalent.

\begin{enumerate}
\item  $\left( sPMV\right)$ holds,

\item  $\left( VC\right) ,\left( TC\right) $ and $\left(
MV\right)  ,$

\item  $\left( VC\right) ,\left( TC\right) $ and the local
diagonal upper estimate
\begin{equation}
P_{n}(x,x)\leq \frac{C\mu (x)}{V(x,e(x,n))},  \label{LDUE}
\end{equation}
holds,

\item  $\left( VC\right) ,\left( TC\right) $ and the local
upper estimate
\begin{equation}
p_{n}\left( x,y\right) \leq \frac{C}{\sqrt{V\left( x,e\left( x,n\right)
\right) V\left( y,e\left( y,n\right) \right) }}\exp \left( -c\kappa
_{3}\left( n,x,y\right) \right) ,    \label{LUE}
\end{equation}
holds.
\end{enumerate}
\end{theorem}

The proof of the diagonal upper estimate can be given along the lines of
\cite{GT2} while the off-diagonal estimate based on a generalization of an
inequality due to Davies \cite{D1}.

\begin{remark}
\bigskip It can also be shown that $\left( p_{0}\right),\left( VC\right) $ and $\left( TC\right) $ imply
\begin{equation}
p_{2n}\left( x,x\right) \geq \frac{c}{V\left( x,e\left( x,2n\right) \right) }%
.
\end{equation}
\bigskip
\end{remark}

\bigskip One gets a weaker upper estimate introducing
\begin{equation}
\kappa \left( n,x,y\right) =\min\limits_{z\in A_{x,y}}\left\{ k_{z}\left( 3n,%
\frac{1}{3}d\left( x,y\right) \right) \right\}
\end{equation}
where $A_{x,y}=B\left( x,d\left( x,y\right) \right) \cup B\left( y,d\left(
x,y\right) \right) $ \ if $d\left( x,y\right) >3\left[ e\left( x,n\right)
+e\left( y,n\right) \right] $ and $\kappa \left( n,x,y\right) =0$ \
otherwise. Similarly we introduce
\begin{equation}
l\left( n,x,y\right) =\max_{z\in A_{x,y}}\left\{ k_{z}\left( n,d\left(
x,y\right) \right) \right\} .
\end{equation}
\ Let us measure the inhomogeneity of the mean exit time for any $A\subset
\Gamma $ by
\begin{equation}
\delta \left( n,A\right) =\log \left[ \max_{z,v\in A}\frac{e\left(
z,n\right) }{e\left( v,n\right) }\right]
\end{equation}
and denote the lower sub-Gaussian kernel by $\nu ;$%
\begin{equation}
\nu \left( n,x,y\right) =l\left( n,x,y\right) \left[ 1+\delta \left(
n,A_{x,y}\right) \right] .
\end{equation}

\begin{definition}
The weighted graph $\left( \Gamma ,\mu \right) $ satisfies the elliptic
Harnack inequality $(\mathbf{H})$ if there is a $C>0$ such that for all $%
x\in \Gamma $ and $R>0$ and for all $u\geq 0$ on $\overline{B}\left(
x,2R\right) $ harmonic functions on $B(x,2R)$ which means that
\begin{equation}
Pu=u
\end{equation}on $B(x,2R)$,
the following inequality holds
\begin{equation}
\max_{B(x,R)}u\leq C\min_{B(x,R)}u.
\end{equation}
\end{definition}

Using the above notations the following statement can be given, which on the
upper estimate side is direct consequence of the above results observing
that the elliptic Harnack inequality implies the mean value inequality $%
\left( MV\right)  .$

\begin{theorem}
\label{tLsGE}Assume that $\left( \Gamma ,\mu \right) $ satisfies
$\left( p_{0}\right)  ,\left( VC\right) ,\left( TC\right) $ and
the elliptic Harnack inequality $\left( H\right) $, then
\begin{equation}
p_{n}\left( x,y\right) \leq \frac{C}{\sqrt{V\left( x,e\left( x,n\right)
\right) V\left( y,e\left( y,n\right) \right) }}\exp \left( -c\kappa \left(
n,x,y\right) \right) ,  \label{UE2}
\end{equation}
and
\begin{equation}
\widetilde{p}_{n}\left( x,y\right) \geq \frac{c}{V\left( x,e\left(
x,n\right) \right) }\exp \left( -C\nu \left( n,x,y\right) \right) .
\label{LE2}
\end{equation}
where $\widetilde{p}_{n}=p_{n+1}+p_{n}.$
\end{theorem}

\begin{remark}
One can rewrite $\left( \ref{LE2}\right) $ in the form of
\begin{equation}
\widetilde{p}_{n}\left( x,y\right) \geq \frac{c}{\sqrt{V\left( x,e\left(
x,n\right) \right) V\left( y,e\left( y,n\right) \right) }}\exp \left( -C\nu
\left( n,x,y\right) \right) ,
\end{equation}
to be compared to $\left( \ref{UE2}\right) .$
\end{remark}

The lower estimate is proved via an important intermediate estimate (called
near diagonal lower estimate c.f. \cite{GT2} or \cite{TD}) \ then the
standard Aronson's chaining argument can be used.

\begin{remark}
\label{rmeet}One should recognize that the upper and lower estimate rely on
comparison of volume and exit times of a chain of balls connecting $x$ \ and
$y.$ \ If the mean exit time is basically independent of the \ center of the
ball it is clear from the definitions that $\kappa \simeq l$ , $\delta
\simeq 1$ \ and hence $\kappa \simeq \nu $ which means that the upper and
lower estimate are the same up to the constants.
\end{remark}

\section{Two-sided estimates}

The semilocal framework is received from the local one if we assume that.
\begin{equation}
E\left( x,R\right) \simeq E\left( y,R\right)  \label{E}
\end{equation}
for all $x,y\in \Gamma .$ \ The study of semi-local situation starts with
the investigation of the space-time scale function $F(R),$ $R\geq 0$ which
is
\begin{equation}
F(R)=\inf_{x\in \Gamma }E(x,R).  \label{EF}
\end{equation}
From $\left( E\right) _{\left( \ref{E}\right) }$ it follows that $F$
satisfies with a fixed $C_{0}>1$ for all $x\in \Gamma $ and $R\geq 0$%
\begin{equation}
F(R)\leq E(x,R)\leq C_{0}F(R).
\end{equation}
Function $F$ inherits certain properties of $E(x,R),$ among others from $%
\left( TD\right)  $ it follows that
\begin{equation}
F(2R)\leq D_{F}F(R).  \label{ED1F}
\end{equation}
The inherited properties are referred by the notation $\left(
ED_{F}\right) $ \ (c.f. \cite{TD}). This function takes over the
role of $R^{\beta }\left( \textnormal{ or }  R^{2}\right) .$ \ \
The inverse function of $F,$ $f(.)=F^{-1}(.)$ takes over the role
of $R^{\frac{1}{\beta }}$ $\left( R^{\frac{1}{2}}\right) $ in the
(sub-)Gaussian estimates.

\begin{definition}
The transition probability satisfies the sub-Gaussian upper estimate $\left(
\mathbf{UE}_{F}\right) $ with respect to $F$ if there are $%
c,C>0$ such that
\begin{equation}
P_{n}(x,y)\leq \frac{C\mu (y)}{V(x,f(n))}\exp -ck(n,d(x,y)),  \label{UEF}
\end{equation}
and the sub-Gaussian lower estimate $\left(
\mathbf{LE}_{F}\right) $ is satisfied if
\begin{equation}
\widetilde{P}_{n}(x,y)\geq \frac{c\mu (y)}{V(x,f(n))}\exp -Ck(n,d(x,y)),
\label{LEF}
\end{equation}
where $\widetilde{P_{n}}=P_{n}+P_{n+1}$ and the kernel function $%
k=k(n,R)\geq 1$, is defined as the maximal integer for which
\begin{equation}
\frac{n}{k}\leq qF(\left\lfloor \frac{R}{k}\right\rfloor )  \label{kdef}
\end{equation}
or $k=0$ by definition if there is no appropriate $k.$
\end{definition}

As we indicated the parabolic and elliptic Harnack inequalities play
important\ role in the study of two-sided bound of the heat kernel. \ Here
we give their formal definitions of the parabolic one..

\begin{definition}
The weighted graph $\left( \Gamma ,\mu \right) $ satisfies the ($F-$%
parabolic or simply$)$ parabolic Harnack inequality $\mathbf{(PH}_{F}\mathbf{%
)}$ if \ the following condition holds. There is a $C_{H}>0$ constant \ such
that for any solution $u\geq 0$ of the equation
\begin{equation}
u_{n+1}(x)=Pu_{n}(x)
\end{equation}
on $\mathcal{U}=[k,k+F(4R)]\times B(x,2R)$ for $k,R\in N$ the
following is true. On the smaller cylinders defined by
\begin{eqnarray*}
\mathcal{U}^{-} &=&[k+F(R),k+F(2R)]\times B(x,R)\textnormal{ } \\
\textnormal{and }\mathcal{U}^{+} &=&[k+F(3R),k+F(4R)]\times B(x,R)
\end{eqnarray*}
and taking $(n_{-},x_{-})\in \mathcal{U}^{-},(n_{+},x_{+})\in \mathcal{U}%
^{+},d(x_{-},x_{+})\leq n_{+}-n_{-}$ the inequality
\begin{equation}
u(n_{-},x_{-})\leq C_{H}\widetilde{u}(n_{+},x_{+})
\end{equation}
holds, where $\widetilde{u}_{n}=u_{n}+u_{n+1}$ short notation was used. The
elliptic Harnack inequality is a direct consequence of the F-parabolic one
as it is true for the classical case.
\end{definition}

Based on the above definitions the following theorem can be formulated.

\begin{theorem}
\label{tmain+}If a weighted graph $(\Gamma ,\mu )$ satisfies
$\left( p_{0}\right)  $ then the following statements are
equivalent.

\begin{enumerate}
\item  $\exists F$ satisfies $(ED_{F})$ and the $F$-parabolic Harnack
inequality $\left( PH_{F}\right) $,

\item  $\exists F$ satisfies $(ED_{F}),$ $(UE_{F})$ and $\left( LE_{F}\right),$

\item  $(VD),\left( TD\right),(E)$ and $\left( H\right) $ hold.
\end{enumerate}
\end{theorem}

Theorem \ref{tLsGE} implies in the semi-local framework
the corresponding off-diagonal estimate ( see also Remark \ref{rmeet}). The
other implication have been proved in \cite{TD}.

The presented results generalize several works, among others Moser \cite{M1}%
, \cite{M2}, Davies \cite{D}, Coulhon, Grigor'yan \cite{CG}, Grigor'yan \cite
{GA}, Li,Yau \cite{LY}, Varopoulos \cite{V}),\cite{De}, Fabes, Stroock \cite
{FS}, Hebisch, Saloff-Coste \cite{HS}. Let us mention that in \cite{HS} the
equivalence of $\left( PH_{F}\right) $ and the $F$-based two-sided
sub-Gaussian estimate was already shown.


\begin{thebibliography}{1}
\bibitem{B} M.T. Barlow, Diffusion on Fractals, Lecture Notes Math. 1690, Springer 1998, 1-121
\bibitem{BCG} M.T. Barlow, T. Coulhon, A. Grigor'yan, Manifolds and graphs with slow heat kernel decay, Invent. Math., 144 (2001) 609-649
\bibitem{Co1} T. Coulhon, Analysis on infinite graphs with regular volume growth, JE 2070, No 17-18, November 1997, Universite de Cergy-Pontoise
\bibitem{CG} T. Coulhon, A. Grigor'yan, Random walks on graphs with regular volume growth, Geometry and Functional Analysis, 8, (1998) 656-701
\bibitem{D} E.B. Davies, Heat kernels and spectral theory, Cambridge University
Press, Cambridge, 1989
\bibitem{D1}  E.B.Davies, Heat kernel bounds, conservation of probability and the Feller property, J.  d'Analyse Math. 58. (1992), 99-119.
\bibitem{De}  T. Delmotte, Parabolic Harnack inequality and estimates of Markov chains on graphs.Revista Matem\'{a}tica Iberoamericana 1,(1999), 181--232.
\bibitem{FS} E. Fabes, D. Stroock, A new proof of the Moser's parabolic Harnack inequality using the old ideas of Nash, Arch. Rat. Mech. Anal., 96, (1986) 327-338
\bibitem{GA} A. Grigor'yan, The heat equation on non-compact Riemannian manifolds, (in Russian)  Matem. Sbornik 182:1, (1991), 55-87 Engl. transl., Math. USSR db. 72:1 (1992) 47-77
\bibitem{Gr1} A. Grigor'yan, Heat kernels on manifolds, graphs and fractals, Prog. in Math., 201 (2001) 393-405
\bibitem{GT1} A. Grigor'yan, A. Telcs, Sub-Gaussian estimates of heat kernels on infinite graphs, Duke Math. J., 109, 3, (2001), 452-510
\bibitem{GT2} A. Grigor'yan, A. Telcs, Harnack inequalities and sub-Gaussian estimates for random walks, Math. Ann. 324. (2002) 521-556
\bibitem{GM} M. Gromov, Groups of polynomial growth and expanding maps. Publ. Math. Inst. H. Poincar\'{e} Probab. Statist. 53 (1981), 57-73
\bibitem{HK} B. Hambly, T. Kumagai, Heat kernel estimates for symmetric random walks ona class of fractal graphs and stability under rough isometries, preprint
\bibitem{HS} W. Hebisch, L. Saloff-Coste, On the relation between elliptic and parabolic Harnack inequalities, Ann. Inst. Fourier 51 (2001) 5, 1437-1481
\bibitem{LY} P. Li, S.-T. Yau, On the parabolic kernel of the Schr\"{o}dinger operator , Acta Math. 156, (1986) 153-201
\bibitem{M1} J. Moser, On Harnack's Theorem for elliptic differential equations, Communications of Pure and Applied Mathematics, 16, (1964) 101-134
\bibitem{M2} J. Moser, On Harnack's theorem for parabolic differential equations, Communications of Pure and Applied Mathematics, 24, (1971) 727-740
\bibitem{TD} A. Telcs, Volume and time doubling of graphs and random walk, the strongly recurrent case, Communication on Pure and Applied Mathematics, LIV, (2001), 975-1018
\bibitem{V} R. Th. Varopoulos,  Hardy-Littlewood theory for semigroups, J. Functional Analysis 63, (1985) 215-239
\bibitem{VCS} R. Th. Varopoulos, L. Saloff-Coste, T. Coulhon, Analysis and geometry on Groups, Cambridge University Press, 1993 
\bibitem{W} W. Woess, Random walks on infinite graphs and groups, Cambridge University Press, Cambridge, 2000
\end{thebibliography}
\end{document}